**Michiel Hazewinkel**　　　　1　　　　CWI
Direct line: +31-20-5924204　　　　POBox 94079
Secretary: +31-20-5924233　　　　1090GB  Amsterdam
Fax: +31-20-5924166
E-mail: mich@cwi.nl　　　　original version: 13 June 2002
　　　　revised version: 19.10.03

# Explicit polynomial generators for the ring of quasisymmetric functions over the integers

by

*Michiel Hazewinkel, CWI Amsterdam*, <mich@cwi.nl>

**Abstract**. In [5, 6] it has been proved that the ring of quasisymmetric functions over the integers is free polynomial, see also [4]. This is a matter that has been of great interest since 1972; for instance because of the role this statement plays in a classification theory for noncommutative formal groups that has been in development since then, see [2] and [9] and the references in the latter. Meanwhile quasisymmetric functions have found many more aplications, [3]. However, the proofs in [5, 6] do not give explicit polynomial generators for *QSymm* over the integers. In this note I give a (really quite simple) set of polynomial generators for *QSymm* over the integers.

**MSCS**: 14L05, 05Exx

**Key words and key phrases**: Quasisymmetric function, symmetric function, plethysm, $\lambda$ – ring, lambda ring, Frobenius operator, Adams operator.

**Générateurs explicites pour les fonctions quasisymmétriques sur les entiers rationels.**

**Resumé.** Dans [5, 6] il a été démontré que l'anneau de fonctions quasisymétriques est polynomialement libre sur l'anneau de base **Z**. C'est là une question importante etudiée depuis 1972; par exemple cet énoncé joue un rôle important dans la théorie de la classification des groupes formel noncommutatifs, voir [2]. [9] et les références données en [9]. Entretemps, les fonctions quasisymétriques ont reçu beaucoup d'applications, [3]. Par contre les démonstrations données dans [5, 6] ne fournissent pas des générateurs polynomiaux explicites pour *QSymm* sur l'anneau des entiers rationels. Dans cette Note nour présentons un ensemble (vraiment très simple) de générateurs polynomiaux pour *QSymm* sur **Z**

**Version française abrégée**

Le but de cette Note est de donner des générateurs explicites pur l'anneau des fonctions quasisymétriques sur l'anneau des entiers rationels **Z**. Alors, soit, **Z**$[x_1, x_2, \cdots]$ l'anneau des polynômes a coefficients entiers en $x_1, x_2, \cdots$ et, comme d'habitude, soient

　　　　*Symm*　　*QSymm*　　**Z**$[x_1, x_2, \cdots]$

les sous-anneaux des fonctions symétriques et quasisymétriques. Si $\alpha = [a_1, \cdots, a_m]$ est une composition de $n$, ça veut dire une suite d'entiers naturels $a_1, a_2, \cdots, a_m$ telle que $\sum_i a_i = n$, la fonction quasissymétrique monomiale associée est

**Michiel Hazewinkel**
Direct line: +31-20-5924204
Secretary: +31-20-5924233
Fax: +31-20-5924166
E-mail: mich@cwi.nl

$$\alpha = \sum_{j_1 < \cdots < j_m} x_{j_1}^{a_1} x_{j_2}^{a_2} \cdots x_{j_m}^{a_m}$$

Ces fonctions quasisymétriques monomiales forment une base de *QSymm* (comme groupe Abélien libre).

Il y a une structure de $\lambda$-anneau sur $\mathbf{Z}[x_1, x_2, \cdots]$, donnée par

$$\lambda_i(x_j) = \begin{cases} x_j & \text{if } i = 1 \\ 0 & \text{if } i = 0 \end{cases}, \quad j = 1, 2, \cdots$$

Les sous-anneaux *Symm* et *QSymm* sommes stable sous ces opérations.

Un mot de Lyndon $\alpha = [a_1, \cdots, a_m]$ est dite élémentaire si $(a_1, a_2, \cdots, a_m) = 1$. Soit *eLYN* l'ensemble des mots de Lyndon élémentaires.

*Théorème.* Les $e_n(\alpha) = \lambda^n(\alpha)$, $\alpha \in eLYN$, $n = 1, 2, \cdots$ forment une base polynomiale libre pour *QSymm* sur $\mathbf{Z}$.

Si $\alpha = [1]$, $e_n(\alpha) = e_n$, la fonction symétrique élémentaire de degrée $n$, ce qui peut expliquer la notation.

Les opérateurs de Adams associés à la structure de $\lambda$-anneau sur *QSymm* sont données sur les fonctions quasisymétriques monomiales par

$$\mathbf{f}_n([a_1, \cdots a_m]) = [na_1, \cdots na_m]$$

Il y a des formules bien connues qui relient les $\lambda_i$ et les $\mathbf{f}_i$, ce qui permet d'écrire des expressions explicites pour les $e_n(\alpha)$.

## 0. Introduction.

As indicated in the abstract a somewhat important problem is the finding of explicit free polynomial generators for the ring of quasisymmetric functions over the integers. A seminal inspirational formula for this was (and is) the following observation

$$\exp(\sum_{n=1}^{\infty} n^{-1}[na_1, na_2, \cdots, na_m]t^n) \in QSymm[[t]] \quad (0.1)$$

It will not be very clear from what follows just what this formula has to do with the actual proof below. (But see remark 2.4.) The fact that this formula is actually of central importance is quite effectively hidden in the use of plethysms as employed in section 2.

A much more detailed paper explaining this and much more is in preparation.

## 1. Lambda ring structure on *QSymm*.

Consider the rings of symmetric functions and quasisymmetric functions in infinitely many variables $x_1, x_2, \cdots$ over the integers

$$Symm = \mathbf{Z}[e_1, e_2, \cdots] \quad QSymm \quad \mathbf{Z}[x_1, x_2, \cdots] \tag{1.1}$$

Here the $e_i$ are the elementary symmetric functions in the $x_j$. For some details on quasisymmetric functions and definitions of various concepts used below see [4, 5, 6]. There is a well known $\lambda$ – ring structure on $\mathbf{Z}[x_1, x_2, \cdots]$ given by

$$\lambda_i(x_j) = \begin{cases} x_j & \text{if } i = 1 \\ 0 & \text{if } i = 0 \end{cases}, \quad j = 1, 2, \cdots \tag{1.2}$$

The asscociated Adams operators, determined by the formula

$$t \frac{d}{dt} \log \lambda_t(a) = \sum_{n=1} (-1)^n \mathbf{f}_n(a) t^n \tag{1.3}$$

where

$$\lambda_t(a) = 1 + \sum_{n=1} \lambda_n(a) t^n \tag{1.4}$$

are the ring endomorphisms

$$\mathbf{f}_n : x_j \mapsto x_j^n \tag{1.5}$$

There are well-known determinantal relations between the $\lambda_n$ and the $\mathbf{f}_n$ as follows

$$n! \lambda_n(a) = \det \begin{pmatrix} \mathbf{f}_1(a) & 1 & 0 & \cdots & 0 \\ \mathbf{f}_2(a) & \mathbf{f}_1(a) & 2 & \ddots & \vdots \\ \vdots & \vdots & \ddots & \ddots & 0 \\ \mathbf{f}_{n-1}(a) & \mathbf{f}_{n-2}(a) & \cdots & \mathbf{f}_1(a) & n-1 \\ \mathbf{f}_n(a) & \mathbf{f}_{n-1}(a) & \cdots & \mathbf{f}_2(a) & \mathbf{f}_1(a) \end{pmatrix} \tag{1.6}$$

$$\mathbf{f}_n(a) = \det \begin{pmatrix} \lambda_1(a) & 1 & 0 & \cdots & 0 \\ \lambda_2(a) & \lambda_1(a) & 1 & \ddots & \vdots \\ \vdots & \vdots & \ddots & \ddots & 0 \\ (n-1)\lambda_{n-1}(a) & \lambda_{n-2}(a) & \cdots & \lambda_1(a) & 1 \\ n\lambda_n(a) & \lambda_{n-1}(a) & \cdots & \lambda_2(a) & \lambda_1(a) \end{pmatrix} \tag{1.7}$$

(which come from the Newton relations between the elementary symmetric functions and the power sum symmetric functions).

It follows that the subrings *Symm* and *QSymm* are stable under the $\lambda_n$ and $\mathbf{f}_n$ because

$\lambda_n(QSymm \otimes_\mathbf{Z} \mathbf{Q}) \subset QSymm \otimes_\mathbf{Z} \mathbf{Q}$ by (1.6) and $(QSymm \otimes_\mathbf{Z} \mathbf{Q}) \cap \mathbf{Z}[x_1, x_2, \cdots] = QSymm$, and similarly for *Symm*. It follows that *QSymm* and *Symm* have induced $\lambda$ – ring structures.

A composition $\alpha = [a_1, a_2, \cdots, a_m]$, $a_i \in \mathbf{N} = \{1, 2, \cdots\}$ defines a monomial quasisymmetric function also denoted by $\alpha$

$$\alpha = \sum_{i_1 < i_2 < \cdots < i_m} x_{i_1}^{a_1} x_{i_2}^{a_2} \cdots x_{i_m}^{a_m} \tag{1.8}$$

The empty composition corresponds to the quasisymmetric function 1. The monomial quasisymmetric functions form a basis for *QSymm* as an Abelian group.

1.9. *Lemma.* $\mathbf{f}_n([a_1, a_2, \cdots, a_m]) = [na_1, na_2, \cdots, na_m]$

This follows immediately from (1.8) and (1.5).

Give a composition $\alpha = [a_1, a_2, \cdots a_m]$ weight $\mathrm{wt}(\alpha) = a_1 + a_2 + \cdots + a_m$ and length $\mathrm{lg}(\alpha) = m$. The wll-ordering on compositions is defined by "weight first, then length, then lexicographic". Thus for instance

$$[5] >_{wll} [1,1,2] >_{wll} [2,2] >_{wll} [1,3]$$

1.10. *Lemma.* If $\alpha$ is a Lyndon word (= Lyndon composition)

$$\lambda_n(\alpha) = \alpha^{*n} + (\text{wll} - \text{smaller than } \alpha^{*n}) \tag{1.11}$$

where $*$ denotes concatenation (of compositions) and (wll – smaller than $\alpha^{*n}$) stands for a $\mathbf{Z}$-linear combination of monomial quasisymmetric functions that are wll-smaller than $\alpha^{*n}$.

This follows immediately from formula (1.6) and lemma 1.9. Indeed, expanding the determinant (1.6) we see that

$$n!\lambda_n(\alpha) = \alpha^{*n} + (\text{monomials of length} \le (n-1) \text{ in the } \mathbf{f}_i(\alpha))$$

All monomials occurring in $\lambda_n(\alpha)$ are of equal weight $n\mathrm{wt}(\alpha)$. By the formula above the longest ones come from the power $\alpha^{*n}$ and are of length $n\,\mathrm{lg}(\alpha)$. Because $\alpha$ is Lyndon the lexicographic largest term of these is $\alpha^{*n}$ and it occurs with coefficient $n!$.

For any $\lambda$ – ring $R$ there is an associated mapping

$$Symm \times R \longrightarrow R, \quad (\varphi, a) \mapsto \varphi(\lambda_1(a), \lambda_2(a), \cdots, \lambda_n(a), \cdots) \tag{1.12}$$

I.e. write $\varphi \in Symm$ as a polynomial in the elementary symmetric functions $e_1, e_2, \cdots$ and then substitute $\lambda_i(a)$ for $e_i$, $i = 1, 2, \cdots$. For fixed $a \in R$ this is obviously a homomorphism of rings $Symm \longrightarrow R$. We shall often simply write $\varphi(a)$ for $\varphi(\lambda_1(a), \lambda_2(a), \cdots, \lambda_n(a), \cdots)$. Another way to see (1.12) is to observe that for fixed $a \in R$ $(\varphi, a) \mapsto \varphi(\lambda_1(a), \lambda_2(a), \cdots) = \varphi(a)$ is the unique homomorphism of $\lambda$-rings that takes $e_1$ into $a$. (*Symm* is the free $\lambda$-ring on one generator, see also [8]. Note that

$$e_n(\alpha) = \lambda_n(\alpha), \quad p_n(\alpha) = \mathbf{f}_n(\alpha) = [na_1, na_2, \cdots, na_m] \tag{1.13}$$

The first formula of (1.13) is by definition and the second follows from (1.7) because the relations between the $e_n$ and $p_n$ are precisely the same as between the $\lambda_n(a)$ and the $\mathbf{f}_n(a)$ (see (1.3)).

 1.14. *Lemma.* For any $\varphi, \psi \in Symm$ and $a \in R$, where $R$ is a $\lambda$-ring

$$\varphi(\psi(a)) = (\varphi \circ \psi)(a) \tag{1.15}$$

where $\varphi \circ \psi$ is the (outer) plethysm of $\varphi, \psi \in Symm$.

This is well known, see [8], p. 134, remark 1. For our purposes here it does not matter just how plethysm is defined. The only thing needed is that there is some element $\varphi \circ \psi \in Symm$ such that (1.15) holds.

## 2. Explicit polynomial generators for *QSymm*.

Let *LYN* denote the set of Lyndon words and let *eLYN* be the set of elementary (or reduced) Lyndon words, i.e. the set of those Lyndon words $\alpha = [a_1, a_2, \cdots, a_m]$ for which $g(\alpha) = \gcd\{a_1, a_2, \cdots, a_m\} = 1$.

 2.1. *Theorem.* The $e_n(\alpha)$, $\alpha \in eLYN$, $n \in \mathbf{N}$ form a free set of polynomial generators over the integers for the ring of quasisymmetric functions *QSymm*.

Proof. The difficult part is to prove generation, i.e. that every basis element $\beta$ (in the Abelian group sense) of *QSymm* can be written as a polynomial in the $e_n(\alpha)$, $\alpha \in eLYN$, $n \in \mathbf{N}$. The rest follows by the same counting argument that was used in, e.g., [5]. So let us prove generation.
 To start with, let $\beta = [b_1, b_2, \cdots, b_m]$ be a Lyndon composition. Then taking $\alpha = \beta_{red} = [g(\beta)^{-1}b_1, g(\beta)^{-1}b_2, \cdots, g(\beta)^{-1}b_m]$, and $n = g(\beta)$ we have, using (1.13), $\beta = p_n(\alpha)$ which is a polynomial in the $e_n(\alpha)$, and thus $\beta \in R$.
 We now proceed by induction on the wll-ordering. The case of weight 1 is trivial. For each separate weight the induction starts because of what has just been said because compositions of length 1 are Lyndon.
 So let $\beta$ be a composition of weight $\geq 2$ and length $\geq 2$. By the Chen-Fox-Lyndon concatenation factorization theorem, [1]

$$\beta = \beta_1^{*r_1} * \beta_2^{*r_2} * \cdots * \beta_k^{*r_k}, \quad \beta_i \in LYN, \quad \beta_1 >_{lex} \beta_2 >_{lex} \cdots >_{lex} \beta_k \tag{2.2}$$

where, as before, the $*$ denotes concatenation and $\beta >_{lex} \beta'$ means that $\beta$ is lexicographically strictly larger than $\beta'$.
 If $k \geq 2$, take $\beta' = \beta_1^{*r_1}$ and for $\beta''$ the corresponding tail of $\beta$ so that $\beta = \beta' * \beta''$.

Then

$$\beta'\beta'' = \beta'\ \beta'' + (\text{wll} - \text{smaller than } \beta) = \alpha + (\text{wll} - \text{smaller than } \beta)$$

and with induction it follows that $\beta \in R$.

There remains the case that $k = 1$ in the CFL-factorization (2.2). In this case take $\alpha = (\beta_1)_{red}$ and observe that by Lemma 1.10 and (1.13)

$$\beta = e_n(p_{g(\beta_1)}(\alpha)) + (\text{wll-smaller than } \beta) \tag{2.3}$$

On the other hand, by Lemma 1.14

$$e_n(p_{g(\beta_1)}(\alpha)) = (e_n \circ p_{g(\beta_1)})(\alpha)$$

Where $e_n \circ p_{g(\beta_1)}$ is some polynomial with integer coefficients in the $e_j$, and hence $(e_n \circ p_{g(\beta_1)})(\alpha)$ is a polynomial with integer coefficients in the $e_j(\alpha)$. With induction this finishes the proof.

2.4. *Remark*. In terms of the $e_n(\alpha)$, for any composition $\alpha = [a_1, a_2, \cdots, a_m]$ the exponential from the introduction is equal to

$$\exp(\sum_{n=1} (-1)^{n-1} n^{-1} [na_1, na_2, \cdots, na_m] t^n) = 1 - e_1(\alpha)t + e_2(\alpha)t^2 - e_3(\alpha)t^3 + \cdots$$

**Acknowledgement.** Define the quasisymmetric functions $a_n(\alpha)$ by the triangular relation

$$\prod_n (1 - a_n(\alpha)t^n) = 1 - e_1(\alpha)t + e_2(\alpha)t^2 - e_3(\alpha)t^3 + \cdots$$

Of course the realations between $a_n(\alpha)$ and the $e_n(\alpha)$ are given by triangular matrices with $\pm 1$ on the diagonal. So if one is a set of gnerators so is the other. In Spring 2002. E J Ditters wrote me a letter stating that the $a_n(\alpha)$ are a free polynomial set of generators of *QSymm*. The proof rests on observing that these elements when evaluated on the explicit basis of Prim(*NSymm*) from [6} give a triangular matrix with $\pm 1$ on the diagonal. Immediately afterwards I gave a proof that did not use duality but still used the abstract theorem that *QSymm* was free. Some months later I found the proof that is in this note. Thus, historicaly speaking, this is the fifth proof that *QSymm* is freely generated and the third one that also provides an explicit set of generators.